***

# Game theory: An appropriate tool for strategic analysis ?

- Didier BAZALGETTE, PhD, Adjunct faculty SciencesPo, Associate Researcerh Center for Doctrine and Command Teaching
- Marc-Olivier BOISSET, PhD, Associate Researcher, Military College Saint Cyr Coëtquidan
- Jean LANGLOIS-BERTHELOT, PhD, Currently pursuing research at Stanford's Freeman Spogli Institute for International Studies on a Defense Innovation Agency/Ecole Polytechnique Fellowship

According to one of its founders, game theory would soon be a little over a hundred years old[1]. A well-known theory in the academic and strategic fields, it is often the subject of criticism and, for many, seems too simplistic to be truly operational for decision-making in complex strategic environments.[2]

In this synthesis article, we analyze the main criticisms made of game theory in strategy and we show to what extent certain aspects of research in game theory have important potential, precisely when the strategic community takes record of the work carried out by academics and operational staff in recent years in emerging research fields.

## "It's trivial, you know" (von Neumann, summer 1949)

---

[1] MORGENSTERN Oskar, « On Some Criticisms of Game Theory », *Research Paper* n° 8, Princeton University, Econometric Research Program, 1964 (https://www.princeton.edu/~erp/ERParchives/archivepdfs/R8.pdf).

[2] CASSIDY John, « The Triumph (and Failure) of John Nash's Game Theory », *The New Yorker*, 27 mai 2015 (https://www.newyorker.com/news/john-cassidy/the-triumph-and-failure-of-john-nashs-game-theory).

Systematically formalized by John von Neumann and Oskar Morgenstern in 1944 in *Theory of Games and Economic Behavior*[3]. Game theory had, however, been the subject of an article by von Neumann in 1928, entitled "*On the Theory of Games of Strategy*"[4], which took up ideas from Borel and Zermelo[5]. Game theory is above all a formal description of strategic interactions between agents.[6]

Mathematician John Nash defended his thesis at Princeton in 1950[7]. This makes it possible to apprehend in a formalized way the potential results of a game/situation with several individuals in a competitive situation and the optimal solution for the mutual gain of the individuals concerned. As such, the Nash equilibrium is a very simple concept, yet it is a major contribution to the dissemination of the work of Morgenstern and von Neumann. When von Neumann indicates that the Nash equilibrium is trivial, it is not so much a criticism since he himself as well as Morgenstern clearly recognize the simple or even simplifying nature of game theory for strategic decision-making in a situation of complexity[8].

**Complex strategic environment and simplified formalization of issues**

---

[3] NEUMANN (VON) John & MORGENSTERN Oskar, *Theory of Games and Economic Behavior*, Princeton University Press, 2007 (1944)

[4] Leonard Robert J., « From Parlor Games to Social Science: von Neumann, Morgenstern, and the Creation of Game Theory 1928-1944 », *Journal of Economic Literature*, 33(2), juin 1995, p. 730-761.

[5] Vorob'ev Nikolai N., « The Present State of The Theory of Games », *Russian Mathematical Surveys*, 25(2), 1970, p. 77-136.

[6] Nash John F., « Non-Cooperative Games (PhD Thesis, 1951) », in Kuhn Harold W. & Nasar Syvia (dir.), *The Essential John Nash*, Princeton University Press, 2007.

[7] Cassidy J., *op. cit.*

[8] Morgenstern O., *Op. cit.*

Almost sixty years ago, in his famous speech in Toulon (1964)[9], before NATO experts, Professor Oskar Morgenstern (Princeton) analyzed the many criticisms made by his university colleagues and military experts with regard to of game theory. The main criticisms focused on the eminently simplistic character of game theory in a strategic approach which is structurally a situation of complexity[10].

The debate on the validity of game theory for strategy is an old debate which has been very widely documented, in particular, because it took place within several leading institutions of strategic thought such as the RAND Corporation, Princeton University and Stanford University[11].

Professor Charles Zorgbibe (Panthéon-Sorbonne University) wrote in his article *"Theory of games and international relations"* (1980), « The basic postulate of the various "scientific" approaches is the rationality of the behavior of the actors. But is the international political game rational? To the theoretician, the "decision maker" will be tempted to retort that the action is essentially pragmatic (…) More generally, any theoretical "reading" of an international event runs the risk of being contested, because it is "imposed" on the situation that she claims to decipher »[12].

---

[9] Morgenstern O., *Op. cit.*
[10] Merle Marcel, « Sur la "problématique" de l'étude des relations internationales en France », *Revue française de science politique*, 33(3), 1983, p. 403-427 (https://www.persee.fr/doc/rfsp_0035-2950_1983_num_33_3_411242).
[11] Gambarelli Gianfranco & Guillermo Owen, « The Coming of Game Theory », *Essays in Cooperative Games*, Springer, Boston, MA, 2004, p. 1-18.
[12] Zorgbibe Charles, « Théorie des jeux et relations internationales », *Le Monde Diplomatique*, octobre 1980 (https://www.monde-diplomatique.fr/1980/10/ZORGBIBE/35814).

Strategic decision-making in the context of a military operation requires understanding a very large number of factors and repeatedly adjusting the information on the evolution of these different factors, taking into account the choices to be made obviously cannot be reduced to the approaches permitted by game theory[13]. In the same way, a game theory approach to cybersecurity, analyzing the actions to be taken from a game theory approach, risks making the understanding of the issues very limited. Of course, a formalization of the stakes in game theory can be of particular interest in negotiations, as Schelling[14] has brilliantly shown, and this, in particular in a "game" with a voluntary limitation of the factors to be analyzed such as in the models used to explain to decision-makers the main axes of decisions to be taken within the framework of the arm races for example.[15]

Does this mean that game theory would only be of real operational use in situations of negotiation and strategic (and economic) transactions? However, some works seem to promise an even more ambitious use of game theory. When this is the case, is it still truly game theory ?

## What "game theories" to understand complex phenomena?

Game theory is, as we have seen, often analyzed as a simplifying conception of the decision-making approach, but it is also understood as a theory which, in its original classical form, has had its day. As Colin Camerer points out: "Many researchers reject a part of game theory that seems to them long outdated" [16].

---

[13] Rubinstein Ariel, « Comments on the Interpretation of Game Theory », *Econometrica: Journal of the Econometric Society*, 59(4), juillet 1991, p. 909-924.
[14] Dixit Avinash, « Thomas Schelling's contributions to game theory », *The Scandinavian Journal of Economics*, 108(2), juin 2006, p. 213-229.
[15] Schelling Thomas C., « The Strategy of Conflict. Prospectus for a Reorientation of Game Theory », *Journal of Conflict Resolution*, 2(3), 1958, p. 203-264.

To consider that game theory is outdated would be to ignore the important advances in technical and logical fields that have allowed a renewal of game theory and, above all, could allow game theory to be a a real operational tool for decision-making in complex situations.

The rise of computer calculation tools has made it possible to take into account more and more factors and more and more occurrences of interactions between these factors. Thus a "game" and the dynamism of this game will be able to approach more and more precisely with structurally complex environments such as a battlefield or a computer environment for example.

Elie Bursztein's research work is the epitome of how the principles of classical game theory are used to develop an approach to understanding a complex and dynamic environment such as a computer network. Bursztein, now head of computer security research at Google, defended his thesis (funded by the Directorate General for Armaments) at ENS Cachan in 2008 on "anticipation games" adapted to computer security. "Anticipation games" are certainly based on game theory, but with an increase in complexity (precisely made possible by computing) such that it is

---

[16] Camerer Colin F., « Does strategy research need game theory? », *Strategic Management Journal*, 12(S2), 1991, p. 137-152.

legitimate to wonder about the fact that it is still game theory...[17]

Computing power allows, to a large extent, to give a practical response to the formalization of strategic environment with a wide variety of factors and their interactions. The last twenty years have also been revolutionized by the more systematic use of Bayesian statistics for game theory.

There are mainly two types of statistics: frequentists and Bayesians[18]. The first calculate the probability of an outcome based on past elements deemed similar. We owe the introduction of Bayesian statistics in game theory to the Hungarian academic John Charles Harsanyi[19][20]. The Bayesian approach to game theory makes it possible to better formalize strategic environments that are generally evolving and dynamic. This approach revolutionizes the use of game theory in decision analysis and in particular in approaches aimed at strategic anticipation and situations of uncertainty. As Jeffrey and Zhu (2021)[21] remind us, Bayesian statistics make it possible to regularly update the various factors to be taken into account in a strategic analysis situation and to change the importance given to such and such a factor. It is, they add, a more optimal approach than the solution offered by "classic" statistics (i.e. frequentialists) which is therefore becoming increasingly used by specialists.

---

[17] Bursztein Elie : « Anticipation games », PhD diss., Ph. D. thesis, École normale supérieure de Cachan, 2008 ; Camerer C.F., *op. cit.*
[18] Maccarone Lee T. & Cole Daniel G., « Bayesian Games for the Cybersecurity of Nuclear Power Plants », *International Journal of Critical Infrastructure Protection*, 37, juillet 2022, 100493.
[19] Langlois-Berthelot Jean, *Appréhender et Assurer les risques Cyber au sein des Organisations : une Approche Systémique et Opérationnelle*, PhD in Applied Mathematics, October 2021 (https://www.theses.fr/2021EHES0096#)
[20] Aspremont (d') Claude & Hammond Peter J., « John C. Harsanyi » in *Conversations on Social Choice and Welfare Theory*, vol. 1, Springer, 2021, p. 37-48.
[21] Pawlick Jeffrey & Zhu Quanyan, Game Theory for Cyber Deception, Springer International Publishing, 2021.

**To conclude**

Game theory has been used since the 1940s in strategic decision-making. It is certainly the subject of many criticisms, but these are mainly due to the fact that it is an attempt to formalize strategic environments that are precisely structurally difficult to formalize.

However, game theory has been profoundly transformed in recent years, even leading some experts to qualify it differently. Indeed, the potentials offered by computer sciences and the use of Bayesian statistics are profoundly transforming the use of game theory for strategic decision-making in complex environments.

With these profound transformations, game theory promises important applications in strategy. These transformations could go with a certain obsolescence of the classical models of game theory for apprehending strategic environments.